\documentclass[11pt]{article}
\usepackage{etex}
\reserveinserts{28}
\usepackage[all]{xy}

\usepackage{amsfonts}
\usepackage{amsthm}
\usepackage{enumerate}
\usepackage{graphicx}
\usepackage{mathrsfs}
\usepackage{bm}
\usepackage{cite}
\usepackage{amssymb,amsmath} 
\usepackage{makecell}
\usepackage{tikz}
\usepackage{pgf}
\usepackage{tikz}
\usetikzlibrary{patterns}
\usepackage{pgffor}
\usepackage{pgfcalendar}
\usepackage{pgfpages}
\usepackage{shuffle,yfonts}
\usepackage{mathtools}
\DeclareFontFamily{U}{shuffle}{}
\DeclareFontShape{U}{shuffle}{m}{n}{ <-8>shuffle7 <8->shuffle10}{}

\newcommand{\ola}{\overleftarrow}
\newcommand{\ora}{\overrightarrow}

\newcommand{\bfk}{{\boldsymbol{\sl{k}}}}

 \allowdisplaybreaks
\usetikzlibrary{arrows,shapes,chains}
 \allowdisplaybreaks

\catcode`!=11
\let\!int\int \def\int{\displaystyle\!int}
\let\!lim\lim \def\lim{\displaystyle\!lim}
\let\!sum\sum \def\sum{\displaystyle\!sum}
\let\!sup\sup \def\sup{\displaystyle\!sup}
\let\!inf\inf \def\inf{\displaystyle\!inf}
\let\!cap\cap \def\cap{\displaystyle\!cap}
\let\!max\max \def\max{\displaystyle\!max}
\let\!min\min \def\min{\displaystyle\!min}
\let\!frac\frac \def\frac{\displaystyle\!frac}
\catcode`!=12

\let\oldsection\section
\renewcommand\section{\setcounter{equation}{0}\oldsection}

\allowdisplaybreaks

\def\R{\mathbb{R}}

\def\N{\mathbb{N}}

\def\ze{\zeta}

\theoremstyle{plain}
\newtheorem{thm}{Theorem}[section]
\newtheorem{lem}[thm]{Lemma}
\newtheorem{cor}[thm]{Corollary}
\newtheorem{con}[thm]{Conjecture}

\theoremstyle{definition}

\newtheorem{re}[thm]{Remark}

\begin{document}
\title{\bf General Mneimneh-type Binomial Sum involving Harmonic Numbers}
\author{
{Ende Pan$^{a,}$\thanks{Email: 13052094150@163.com}\quad{and}\quad Ce Xu$^{b,}$\thanks{Email: cexu2020@ahnu.edu.cn}}\\[1mm]
a. \small College of Teacher Education, Quzhou University, \\ \small Quzhou 324022, P.R. China\\
b. \small School of Mathematics and Statistics, Anhui Normal University,\\ \small Wuhu 241002, P.R. China
}

\date{}
\maketitle

\noindent{\bf Abstract.} Recently, Mneimneh proved the remarkable identity
\begin{align*}
\sum_{k=0}^n H_k\binom{n}{k} p^k(1-p)^{n-k}=\sum_{i=1}^n \frac{1-(1-p)^i}{i}\quad (p\in [0,1])
\end{align*}
as the main result of a 2023 \emph{Discrete Mathematics} paper, where $H_k:=\sum\nolimits_{i=1}^k 1/i$ is the classical $k$-th harmonic number. Thereafter, Campbell provided several other proofs of Mneimneh's formula as above in a note published in \emph{Discrete Mathematics} in 2023. Moreover, Campbell also considered how Mneimneh's identity may be proved and generalized using the \emph{Mathematica package Sigma}. In particular, he found the generalized Mneimneh's identity
\begin{align*}
\sum_{k=0}^n x^k y^{n-k} \binom{n}{k}H_k =(x+y)^n \left(H_n-\sum_{i=1}^n \frac{y^i (x+y)^{-i}}{i}\right).
\end{align*}
In this paper, we will prove a more generalization of Mneimneh's identity involving Bell numbers and some Mneimneh-type identities involving (alternating) harmonic numbers by using a few results of our previous papers.

\medskip

\medskip


\section{Introduction}
Let $\N=\{1,2,3,\ldots\}$ be the set of natural numbers, and $\N_0:=\N\cup \{0\}$. In his recent paper \cite{M2023}, Mneimneh used the method of probabilistic analysis to establish the following binomial sum identity involving harmonic numbers
\begin{align}\label{Mneimneh-identity}
\sum_{k=0}^n H_k\binom{n}{k} p^k(1-p)^{n-k}=\sum_{i=1}^n \frac{1-(1-p)^i}{i}\quad (p\in [0,1]),
\end{align}
where $H_k$ is the \emph{classical $k$-th harmonic number} defined by
\[H_k:=\sum_{i=1}^k \frac1{i}\quad \text{and}\quad H_0:=0.\]
Quite recently, Campbell \cite{C2023} gave two new proofs of \eqref{Mneimneh-identity} by using Zeilberger's algorithm and beta-type integral formula. Further, using the Sigma package for the Mathematica Computer Algebra System, Campbell gave the more general Mneimneh's identity
\begin{align}\label{MG-Mneimneh-identity}
\sum_{k=0}^n x^k y^{n-k} \binom{n}{k}H_k =(x+y)^n \left(H_n-\sum_{i=1}^n \frac{y^i (x+y)^{-i}}{i}\right).
\end{align}
Clearly, setting $(x,y)=(p,1-p)$ in \eqref{MG-Mneimneh-identity} gives \eqref{Mneimneh-identity}. Moreover,
Campbell also emphasized that this approach may also be used to derive identities for expressions such as
\begin{align*}
\sum_{k=0}^n x^k y^{n-k} \binom{n}{k}H_k^2.
\end{align*}
In a 2024 \emph{Discrete Mathematics} article \cite{KW2024}, Komatsu and Wang extended Mneimneh's formula to the generalized hyperharmonic numbers. Gen$\check{\rm c}$ev \cite{G2024} studied of the binomial sum
\begin{align*}
\sum_{k=0}^n \binom{n}{k}\left(\sum_{j=1}^k \frac{z^j}{j}\right)p^k(1-p)^k\quad (r\in \N,\ p,z\in \R,\ p\neq 1)
\end{align*}
and established explicit formula, see \cite[Thm. 2.1]{G2024}.

In this paper, using the method of integrals of natural logarithms, we will establish the explicit formulas of the following general Mneimneh-type binomial sums involving Bell numbers and some Mneimneh-type binomial sums involving (alternating) harmonic numbers.
\begin{thm}\label{thm:Mneimneh-type} For any reals $x,y$ with $x/(x+y)\geq 0$ and $n,p\in \N$, we have
\begin{align}\label{equ:Mneimneh-type}
\sum_{k=0}^n x^ky^{n-k}\binom{n}{k} Y_p(k)=(-1)^{p-1}p! (x+y)^n \sum_{j=1}^n \left(1-\Big(\frac{y}{x+y}\Big)^j\right)\sum_{i=0}^{p-1} (-1)^i \frac{Y_i(n)}{i!} \frac{s(j,p-i)}{j!},
\end{align}
where $s(n,k)$ and $Y_k(n)$ stand for the (unsigned) Stirling numbers of the first kind and Bell numbers (see Section \ref{Intr-STT-Bell}), respectively.
\end{thm}

\begin{re}
All Stirling numbers $s(n,k)$ and Bell numbers $Y_k(n)$ can be expressed in terms of a linear combination of products of harmonic numbers.
\end{re}

In particular, letting $p=1,2$ in Theorem \ref{thm:Mneimneh-type} and noting the facts that $s( {n,1} ) = ( {n - 1} )!,s( {n,2} ) = ( {n - 1} )!{H_{n - 1}}, Y_0(n)=1,{Y_1}( n ) = {H_n}, {Y_2}( n ) = H_n^2 + {H^{(2)} _n}$, we also obtain \eqref{MG-Mneimneh-identity} and
\begin{align}\label{TMG-Mneimneh-identity}
\sum_{k=0}^n x^k y^{n-k} \binom{n}{k}\big(H_k^2+H_k^{(2)} \big) =2(x+y)^n \sum_{j=1}^n \left(1-\Big(\frac{y}{x+y}\Big)^j\right)\frac{H_n-H_{j-1}}{j}.
\end{align}

\begin{thm} \label{thm:Mneimneh-type-G} For any reals $x,y$ and $z\in (-\infty,1)$ with $x/(x+y)\geq 0$ and $n\in \N$, we have
\begin{align}\label{equ:Mneimneh-type-G}
\sum_{k=0}^n x^ky^{n-k}\binom{n}{k} \left(\sum_{j=1}^k \frac{z^j}{j}\right)=\sum_{j=1}^n \frac{(y+zx)^j-y^j}{j}(x+y)^{n-j}.
\end{align}
\end{thm}

Obviously, letting $z\rightarrow 1$ in \eqref{equ:Mneimneh-type-G} gives \eqref{MG-Mneimneh-identity}. Setting $z=-1$ in Theorem \ref{thm:Mneimneh-type-G} yields the following corollary.
\begin{cor} \label{cor:Mneimneh-type-AHN} For any reals $x,y$ with $x/(x+y)\geq 0$ and $n\in \N$, we have
\begin{align}\label{equ:Mneimneh-type-AHN}
\sum_{k=0}^n x^ky^{n-k}\binom{n}{k} {\bar H}_k=\sum_{j=1}^n \frac{y^j-(y-x)^j}{j}(x+y)^{n-j},
\end{align}
where ${\bar H}_k$ is the \emph{alternating $k$-th harmonic number} defined by
\begin{align}
{\bar H}_k:=\sum_{j=1}^k \frac{(-1)^{j-1}}{j}\quad\text{and}\quad {\bar H}_0:=0.
\end{align}
\end{cor}

Specially, setting $(x,y)=(p,1-p)\ (p\in[0,1])$ in \eqref{equ:Mneimneh-type-AHN} yields the following Mneimneh-type identity
\begin{align}\label{equ:Mneimneh-type-AHN-case}
\sum_{k=0}^n {\bar H}_k \binom{n}{k} p^k(1-p)^{n-k}=\sum_{j=1}^n \frac{(1-p)^j-(1-2p)^j}{j}.
\end{align}

From equation \eqref{TMG-Mneimneh-identity} and Theorem \ref{thm:Mneimneh-type-G}, we can also obtain the following corollary.
\begin{cor} \label{cor:Mneimneh-type-order2} For any reals $x,y$ with $x/(x+y)\geq 0$ and $n\in \N$, we have
\begin{align}
&\sum_{k=0}^n x^k y^{n-k} \binom{n}{k}H_k^{(2)} =(x+y)^n \left\{\sum_{j=1}^n \frac{y^j(x+y)^{-j}}{j}\sum_{i=1}^j \frac{y^{-i}(x+y)^i}{i}-\sum_{j=1}^n \frac{y^j(x+y)^{-j}}{j}H_j\right\},\label{TMG-Mneimneh-identity-one}\\
&\sum_{k=0}^n x^k y^{n-k} \binom{n}{k}H^2_k=(x+y)^n \left\{H_n^2+H_n^{(2)}-2\sum_{j=1}^n \frac{y^j(x+y)^{-j}}{j^2}-2H_n \sum_{j=1}^n \frac{y^j(x+y)^{-j}}{j}\atop+3\sum_{j=1}^n \frac{y^j(x+y)^{-j}}{j}H_j-\sum_{j=1}^n \frac{y^j(x+y)^{-j}}{j}\sum_{i=1}^j \frac{y^{-i}(x+y)^i}{i}\right\}.\label{TMG-Mneimneh-identity-two}
\end{align}
\end{cor}
We will prove Theorem \ref{thm:Mneimneh-type} in Section \ref{sec-3}, and Theorem \ref{thm:Mneimneh-type-G} and Corollary \ref{cor:Mneimneh-type-order2} in Section \ref{sec-4}.
\begin{con} The Theorems \ref{thm:Mneimneh-type} and \ref{thm:Mneimneh-type-G} hold for any reals $x,y$ and $z$.

\end{con}

\section{Preliminaries and lemmas}\label{Intr-STT-Bell}

\subsection{(unsigned) Stirling number of the first kind}\label{Intr-Stirling number}
We recall the definition of \emph{(unsigned) Stirling number of the first kind}. Let $s(n,k)$ denote the (unsigned) Stirling number of the first kind, which is defined by \cite{CG1996,C1974}
\begin{align}\label{SN-Dgf}
n!x\left( {1 + x} \right)\left( {1 + \frac{x}{2}} \right) \cdots \left( {1 + \frac{x}{n}} \right) = \sum\limits_{k = 0}^n s(n+1,k+1) x^{k + 1},
\end{align}
with $s(n,k):=0$ if $n<k$, and $s(n,0)=s(0,k):=0,\ s(0,0):=1$, or equivalently, by the generating function:
\begin{align*}
{\log ^k}( {1 - x} ) = {\left( { - 1} \right)^k}k!\sum\limits_{n = 1}^\infty  {s(n,k)\frac{{{x^n}}}{{n!}}}\quad (x \in [- 1,1)).
\end{align*}
The Stirling numbers ${s\left( {n,k} \right)}$ of the first kind satisfy a recurrence relation in the form
\[s( {n,k}) = s({n - 1,k - 1}) + \left( {n - 1} \right)s( {n - 1,k})\quad (n,k \in \N).\]
Obviously, ${s\left( {n,k} \right)}$ can be expressed in terms of a linear combinations of products of harmonic numbers (see \eqref{defn-gener-harmonicnumber}). In particular,
\begin{align*}
& s( {n,1} ) = \left( {n - 1} \right)!,\\& s\left( {n,2} \right) = \left( {n - 1} \right)!{H_{n - 1}},\\& s( {n,3} ) = \frac{{\left( {n - 1} \right)!}}{2}\left[ {H_{n - 1}^2 - {H^{(2)} _{n - 1}}} \right],\\
&s( {n,4} ) = \frac{{\left( {n - 1} \right)!}}{6}\left[ {H_{n - 1}^3 - 3{H_{n - 1}}{H^{(2)} _{n - 1}} + 2{H^{(3)} _{n - 1}}} \right], \\
&s( {n,5} ) = \frac{{\left( {n - 1} \right)!}}{{24}}\left[ {H_{n - 1}^4 - 6{H^{(4)} _{n - 1}} - 6H_{n - 1}^2{H^{(2)} _{n - 1}} + 3(H^{(2)} _{n - 1})^2 + 8H_{n - 1}^{}{H^{(3)} _{n - 1}}} \right].
\end{align*}

In \cite[Thm 2.5]{Xu2017}, we proved
\begin{align}\label{equ-Str-mhs-1}
s\left( {n,k} \right) = \left( {n - 1} \right)!{\zeta _{n - 1}}( {{{\left\{ 1 \right\}}_{k - 1}}})\quad (k,n\in \N),
\end{align}
where $\zeta_n(\{1\}_k)$ is a special multiple harmonic sum defined by ($\{l\}_r$ denotes the sequence obtained by repeating $l$ exactly $r$ times)
\begin{align}
\zeta_n(\{1\}_r):=\sum_{n\geq n_1>n_2>\cdots>n_r>0} \frac{1}{n_1n_2\cdots n_r}.
\end{align}
For $\bfk=(k_1,\ldots,k_r)\in \N^r$ and positive integer $n$, the generalized \emph{multiple harmonic sums} (MHSs) are defined by
\begin{align}
\zeta_n(\bfk)\equiv \zeta_n(k_1,\ldots,k_r):=\sum\limits_{n\geq n_1>\cdots>n_r>0 } \frac{1}{n_1^{k_1}\cdots n_r^{k_r}}\label{MHSs}.
\end{align}
We set ${\zeta _n}(\emptyset ):=1$ if $\bfk=\emptyset$ and ${\zeta_n}(\bfk):=0$ if $n<r$. For $\bfk=(k)\in \N$,
\begin{align}\label{defn-gener-harmonicnumber}
\ze_n(k)\equiv H_n^{(k)}=\sum_{j=1}^n \frac{1}{j^k}
\end{align}
is the $n$-th generalized \emph{harmonic number} of order $k$, and furthermore, if $k=1$ then $H_n\equiv H_n^{(1)}$ is the classical $n$-th harmonic number.
When taking the limit $n\rightarrow \infty$ in \eqref{MHSs} we get the so-called classical \emph{multiple zeta values} (MZVs) (see \cite{H1992,DZ1994}),
\begin{align*}
{\zeta}( \bfk):=\lim_{n\rightarrow \infty}{\zeta_n}(\bfk),
\end{align*}
defined for $\bfk\in \N^r$ and $k_1>1$ to ensure convergence of the series.

\subsection{Bell polynomials}\label{Intr-Bell polynomials}

Define the \emph{exponential partial Bell polynomials} $B_{n,k}$ by
\[
\frac{1}{k!}\left(\sum_{n=1}^\infty x_n\frac{t^n}{n!}\right)^k
    =\sum_{n=k}^\infty B_{n,k}(x_1,x_2,\ldots,x_n)\frac{t^n}{n!}\,,\quad k=0,1,2,\ldots\,,
\]
and the \emph{exponential complete Bell polynomials} $Y_n$ by
\[
Y_n(x_1,x_2,\ldots,x_n):=\sum_{k=0}^nB_{n,k}(x_1,x_2,\ldots,x_n)
\]
(see \cite[Section 3.3]{C1974}). According to \cite[Eq. (2.44)]{Riordan58}, the complete Bell polynomials $Y_n$ satisfy the recurrence
\[
Y_0=1\,,\quad
Y_n(x_1,x_2,\ldots,x_n)=\sum_{j=0}^{n-1}\binom{n-1}{j}x_{n-j}Y_j(x_1,x_2,\ldots,x_j)
    \,,\quad n\geq1\,,
\]
from which, the first few polynomials can be obtained immediately:
\begin{align*}
&Y_0=1\,,\quad
    Y_1(x_1)=x_1\,,\quad
    Y_2(x_1,x_2)=x^2_1+x_2\,,\quad
    Y_3(x_1,x_2,x_3)=x^3_1+3x_1x_2+x_3\,,\\
&Y_4(x_1,x_2,x_3,x_4)=x^4_1+6x^2_1x_2+4x_1x_3+3x^2_2+x_4\,.
\end{align*}

Define the \emph{Bell number} ${Y_k}( n )$ by
\begin{align}
{Y_k}( n ):= {Y_k}( {{H_n},1!{H_n^{(2)}},2!{H_n^{(3)}}, \cdots ,( {k - 1})!{H_n^{(k)}}} ).
\end{align}
Clearly, the Bell number ${Y_k}( n )$ is a rational linear combination of products of harmonic numbers. We have
\begin{align*}
&{Y_1}( n ) = {H_n},\\&{Y_2}( n ) = H_n^2 + {H^{(2)} _n},\\&{Y_3}( n ) =  H_n^3+ 3{H_n}{H^{(2)} _n}+ 2{H^{(3)} _n},\\
&{Y_4}( n ) = H_n^4 + 8{H_n}{H^{(3)} _n} + 6H_n^2{H^{(2)} _n} + 3(H^{(2)} _n)^2 + 6{H^{(4)} _n},\\
&{Y_5}( n ) = H_n^5 + 10H_n^3{H^{(2)} _n} + 20H_n^2{H^{(3)}_n} + 15{H_n}({H^{(2)}_n})^2 + 30{H_n}{H^{(4)} _n}+ 20{H^{(2)} _n}{H^{(3)} _n} + 24{H^{(5)} _n}.
\end{align*}

In \cite[Eq. (2.9)]{Xu2017}, we showed
\begin{align}\label{equ-Str-mhss-1}
\zeta _n^ \star \left( {{{\{ 1\} }_m}} \right) = \frac{1}{{m!}}{Y_m}\left( n \right)\quad (n,m\in \N_0),
\end{align}
where $\zeta_n^\star(\{1\}_r)$ is a special multiple harmonic star sum defined by
\begin{align*}
\zeta_n^\star(\{1\}_r):=\sum_{n\geq n_1\geq n_2\geq \cdots\geq n_r>0} \frac{1}{n_1n_2\cdots n_r}.
\end{align*}
For $\bfk=(k_1,\ldots,k_r)\in \N^r$ and positive integer $n$, the generalized \emph{multiple harmonic star sums} (MHSSs) are defined by
\begin{align}
\zeta^\star_n(\bfk)\equiv \zeta^\star_n(k_1,\ldots,k_r):=\sum\limits_{n\geq n_1\geq \cdots\geq n_r>0 } \frac{1}{n_1^{k_1}\cdots n_r^{k_r}}\label{MHSSs}.
\end{align}
Similarly, we set ${\zeta^\star_n}(\emptyset ):=1$ if $\bfk=\emptyset$ .

\subsection{Some Lemmas}

Next, we present some lemmas, which are useful in the development of our main theorems.

\begin{lem}(\cite[Thm. 2.9]{XYS2016} )\label{lem-Bell-Log} For $n\in \N$ and $p\in \N_0$, we have
\begin{align}
\int\limits_0^1 {{t^{n - 1}}{{\log}^p}(1 - t)} dt = {( { - 1} )^p}\frac{{{Y_p}\left( n \right)}}{n}.
\end{align}
\end{lem}

\begin{lem}(\cite[Thm. 2.2]{Xu2017}) \label{lem-x-mhs-Log} For $n,m\in \N$ and $x\in(-\infty,1)$, we have
\begin{align}\label{equ-log-x-m}
&\int\limits_0^x {{t^{n - 1}}{{\log}^m}\left( {1 - t} \right)} dt = m!\frac{{{{\left( { - 1} \right)}^m}}}{n}\zeta _n^ \star ( {{{\{ 1\} }_m};x} )\nonumber\\
&\quad\quad\quad\quad\quad\quad\quad\quad\quad+\frac{1}{n}\sum\limits_{j =0}^{m - 1} {{{\left( { - 1} \right)}^{j}}j!\binom{m}{j}{{\log}^{m - j}}\left( {1 - x} \right)} \left( {\zeta _n^ \star( {{{\{ 1\} }_j};x}) - \zeta _n^ \star( {{{\{ 1\} }_j}} )} \right),
\end{align}
where $\zeta_n^\star(\{1\}_r)\ (r\in \N_0)$ is defined by
\begin{align*}
\zeta_n^\star(\{1\}_r;x):=\sum_{n\geq n_1\geq n_2\geq \cdots\geq n_r>0} \frac{x^{n_r}}{n_1n_2\cdots n_r}\quad \text{and}\quad \zeta_n^\star(\emptyset;x):=x^n.
\end{align*}
(Note that in \cite[Thm. 2.2]{Xu2017}, the range of values for $x$ is $x\in[-1,1)$, but in fact, the above equation also holds for $x\in (-\infty,-1]$.)
\end{lem}

\begin{lem}(\cite[Thm. 4.3]{XuZhao2020d})\label{lem-PMPLs2}
For ${\bfk}=(k_1,\ldots,k_r)\in \N^r$ and $l\in\N_0,n\in\N$, we have
\begin{align}\label{eq-PMPLs2}
&\sum_{n\geq n_1\geq n_2\geq \cdots\geq n_r>0}  \prod_{j=1}^r \frac{1}{(n_{j}+l)^{k_{j}}}
=(-1)^r \sum_{j=0}^r (-1)^j\ze^\star_{n+l}(\ora\bfk_{\hskip-2pt j}) \ze_l(\ola\bfk_{\hskip-2pt j+1}),
\end{align}
where $\ora\bfk_{\hskip-2pt j}:=(k_1,\ldots,k_j)$ and $\ola\bfk_{\hskip-2pt j}:=(k_r,\ldots,k_j)$ for all $1\le j\le r$.
\end{lem}

\section{Proof of Theorem \ref{thm:Mneimneh-type}}\label{sec-3}
For $x/(x+y)=0$, namely $x=0$, \eqref{equ:Mneimneh-type} is obviously holds. For $x/(x+y)>0$,
applying Lemma \ref{lem-Bell-Log}, the left hand side of \eqref{equ:Mneimneh-type} can be rewritten as
\begin{align}\label{eq-one-proc-main-proof}
\sum_{k=0}^n x^ky^{n-k}\binom{n}{k} Y_p(k)&=(-1)^p \sum_{k=1}^n x^k y^{n-k}\binom{n}{k} k \int_0^1 t^{k-1}\log^p(1-t)dt\nonumber\\
&=(-1)^p\int_0^1 \log^p(1-t) \sum_{k=1}^n x^k y^{n-k}\binom{n}{k} k t^{k-1}dt \nonumber\\
&=(-1)^p\int_0^1 \log^p(1-t) \sum_{k=0}^{n-1} x^{k+1} y^{n-k-1}\binom{n}{k+1} (k+1) t^{k}dt \nonumber\\
&=(-1)^p\int_0^1 \log^p(1-t) \sum_{k=0}^{n-1} x^{k+1} y^{n-k-1}n\binom{n-1}{k} t^{k}dt \nonumber\\
&=(-1)^pnx\int_0^1 \log^p(1-t) \sum_{k=0}^{n-1} (xt)^{k} y^{n-k-1}\binom{n-1}{k}dt \nonumber\\
&=(-1)^pnx\int_0^1 \log^p(1-t)(tx+y)^{n-1}dt\quad\quad\quad(\text{letting}\ u=tx+y)\nonumber\\
&=(-1)^pn\int_y^{x+y} u^{n-1}\log^p\Big(\frac{x+y-u}{x}\Big)du\nonumber\\
&=(-1)^pn\int_y^{x+y} u^{n-1} \left\{\log\Big(\frac{x+y}{x}\Big)+\log\Big(1-\frac{u}{x+y}\Big)\right\}^pdu \nonumber\\
&=(-1)^pn \sum_{l=0}^p \binom{p}{l}\log^{p-l}\Big(\frac{x+y}{x}\Big) \int_y^{x+y} u^{n-1} \log^{l}\Big(1-\frac{u}{x+y}\Big)du\nonumber\\
&=(-1)^pn \sum_{l=1}^p \binom{p}{l}\log^{p-l}\Big(\frac{x+y}{x}\Big) \int_y^{x+y} u^{n-1} \log^{l}\Big(1-\frac{u}{x+y}\Big)du\nonumber\\
&\quad+(-1)^p \Big((x+y)^n-y^n\Big) \log^{p}\Big(\frac{x+y}{x}\Big)\quad\quad\quad(\text{letting}\ u=v(x+y))\nonumber\\
&=(-1)^pn (x+y)^n\sum_{l=1}^p \binom{p}{l}\log^{p-l}\Big(\frac{x+y}{x}\Big) \int_{y(x+y)^{-1}}^{1} v^{n-1} \log^{l}\Big(1-v\Big)dv\nonumber\\
&\quad+(-1)^p \Big((x+y)^n-y^n\Big) \log^{p}\Big(\frac{x+y}{x}\Big).
\end{align}
(In order to utilize Lemmas \ref{lem-Bell-Log} and \ref{lem-x-mhs-Log}, the integration in the second to last row of the above equation needs to satisfy $y(x+y)^{-1}< 1$, namely $x/(x+y)>0$)
Using Lemmas \ref{lem-Bell-Log} and \ref{lem-x-mhs-Log}, by an elementary calculation, the \eqref{eq-one-proc-main-proof} is equal to
\begin{align}\label{eq-two-proc-main-proof}
\sum_{k=0}^n x^ky^{n-k}\binom{n}{k} Y_p(k)&=(x+y)^n \sum_{1\leq i \leq l\leq p} (-1)^{l-i}i!\binom{l}{i}\binom{p}{l} \log^{p-i} \Big(\frac{x}{x+y}\Big)\nonumber\\&\quad\quad\quad\quad\quad\quad\times\left( {\zeta _n^ \star( {{{\{ 1\} }_i}}) - \zeta _n^ \star\Big( {{{\{ 1\} }_i}};y(x+y)^{-1} \Big)} \right)\nonumber\\
&=(x+y)^n \sum_{1\leq i \leq l\leq p} (-1)^{l-i}i!\binom{l}{i}\binom{p}{l} \log^{p-i} \Big(\frac{x}{x+y}\Big)\nonumber\\&\quad\quad\quad\quad\quad\quad\times \sum_{n\geq k_1 \geq \cdots \geq k_i\geq 1} \frac{1-\Big(\frac{y}{x+y}\Big)^{k_i}}{k_1\cdots k_i}.
\end{align}
Noting the fact that the \eqref{eq-two-proc-main-proof} can be rewritten as
\begin{align}\label{eq-three-proc-main-proof}
\sum_{k=0}^n x^ky^{n-k}\binom{n}{k} Y_p(k)&=(x+y)^n \sum_{i=1}^p  (-1)^{i}i!\log^{p-i} \Big(\frac{x}{x+y}\Big) \sum_{n\geq k_1 \geq \cdots \geq k_i\geq 1} \frac{1-\Big(\frac{y}{x+y}\Big)^{k_i}}{k_1\cdots k_i}\nonumber\\&\quad\quad\quad\quad\quad\quad\times \sum_{l=i}^p (-1)^l \binom{l}{i} \binom{p}{l}
\end{align}
and
\begin{align*}
\sum_{l=i}^p (-1)^l \binom{l}{i} \binom{p}{l}&=\left\{ {\begin{array}{*{20}{c}} (-1)^p,
   {\ \ (p=i),} \;\;\; \ \\
  \quad\quad0,\;\;\ \ \ \ {(\text{otherwise}),}  \\
\end{array} } \right.
\end{align*}
we obtain
\begin{align}\label{eq-four-proc-main-proof}
\sum_{k=0}^n x^ky^{n-k}\binom{n}{k} Y_p(k)=p! (x+y)^n \sum_{n\geq k_1 \geq \cdots \geq k_p\geq 1} \frac{1-\Big(\frac{y}{x+y}\Big)^{k_p}}{k_1\cdots k_p}.
\end{align}

In Lemma \ref{lem-PMPLs2}, replacing $n$ by $n-j+1$ and letting $r=p-1$, $(k_1,\ldots,k_r)=(\{1\}_{p-1})$ and $l=j-1$, we easily obtain
\begin{align}
&\sum_{n-j+1\geq i_1\geq \cdots\geq i_{p-1}\geq 1} \frac{1}{(i_1+j-1)\cdots (i_{p-1}+j-1)}\nonumber\\
&\quad\quad\quad=(-1)^{p-1}\sum_{i=0}^{p-1} (-1)^i\ze^\star_n(\{1\}_i)\ze_{j-1}(\{1\}_{p-1-i}).
\end{align}
Hence,
\begin{align}\label{eq-five-proc-main-proof}
\sum_{n\geq k_1 \geq \cdots \geq k_p\geq 1} \frac{1-\Big(\frac{y}{x+y}\Big)^{k_p}}{k_1\cdots k_p}&=\sum_{j=1}^n \frac{1-\Big(\frac{y}{x+y}\Big)^{j}}{j} \sum_{n\geq k_1\geq \cdots\geq k_{p-1}\geq j} \frac{1}{k_1\cdots k_{p-1}}\nonumber\\
&=\sum_{j=1}^n \frac{1-\Big(\frac{y}{x+y}\Big)^{j}}{j}  \sum_{n-j+1\geq i_1\geq \cdots\geq i_{p-1}\geq 1} \frac{1}{(i_1+j-1)\cdots (i_{p-1}+j-1)}\nonumber\\
&=(-1)^{p-1} \sum_{j=1}^n \frac{1-\Big(\frac{y}{x+y}\Big)^{j}}{j} \sum_{i=0}^{p-1} (-1)^i\ze^\star_n(\{1\}_i)\ze_{j-1}(\{1\}_{p-1-i})\nonumber\\
&=(-1)^{p-1} \sum_{j=1}^n \left(1-\Big(\frac{y}{x+y}\Big)^j\right)\sum_{i=0}^{p-1} (-1)^i \frac{Y_i(n)}{i!} \frac{s(j,p-i)}{j!},
\end{align}
where we used the equations \eqref{equ-Str-mhs-1} and \eqref{equ-Str-mhss-1}.

Finally, substituting \eqref{eq-five-proc-main-proof} into \eqref{eq-four-proc-main-proof} yields the desired evaluation \eqref{equ:Mneimneh-type}. Thus, this concludes the proof of Theorem \ref{thm:Mneimneh-type}.

\section{Proofs of Theorem \ref{thm:Mneimneh-type-G} and Corollary \ref{cor:Mneimneh-type-order2}}\label{sec-4}

For $x/(x+y)=0$, namely $x=0$, Theorem \ref{thm:Mneimneh-type-G} and Corollary \ref{cor:Mneimneh-type-order2} are obviously hold. For $x/(x+y)>0$, in \eqref{equ-log-x-m}, setting $m=1$ and replacing $x$ by $z$ gives
\begin{align}\label{equ-log-m-1-x}
\int_0^z {{t^{n - 1}}\log( {1 - t})} dt = \frac{1}{n}\left\{ {{x^n}\log ( {1 - z} ) - \sum\limits_{j = 1}^n {\frac{{{z^j}}}{j}}  - \log ( {1 -z} )} \right\}.
\end{align}
Hence,
\begin{align*}
\sum\limits_{j = 1}^n \frac{z^j}{j}=(z^n-1)\log(1-z)-n\int_0^z t^{n-1}\log(1-t)dt.
\end{align*}

Hence, by a similar argument as in the proof of \eqref{eq-one-proc-main-proof} gives
\begin{align}\label{equ:Mneimneh-type-AHN-case}
\sum_{k=0}^n x^ky^{n-k}\binom{n}{k} \left(\sum_{j=1}^k \frac{z^j}{j}\right)&=-nx \int_0^z \log(1-t)(xt+y)^{n-1}dt+\log(1-z) \left((xz+y)^n -(x+y)^n\right)\nonumber\\
&=-n\int_{y}^{xz+y} u^{n-1}\log\left(\frac{x+y-u}{x}\right)du+\log(1-z) \left((xz+y)^n -(x+y)^n\right)\nonumber\\
&=-n\int_{y}^{xz+y} u^{n-1}\left\{\log\left(\frac{x+y}{x}\right)+\log\left(1-\frac{u}{x+y}\right)\right\}du\nonumber\\
&\quad+\log(1-z) \left((xz+y)^n -(x+y)^n\right)\nonumber\\
&=-\log\left(\frac{x+y}{x}\right)\left((xz+y)^n-y^n\right)+\log(1-z) \left((xz+y)^n -(x+y)^n\right)\nonumber\\
&\quad-n\int_{y}^{xz+y} u^{n-1} \log\left(1-\frac{u}{x+y}\right)du\nonumber\\
&=-\log\left(\frac{x+y}{x}\right)\left((xz+y)^n-y^n\right)+\log(1-z) \left((xz+y)^n -(x+y)^n\right)\nonumber\\
&\quad-n(x+y)^n \int_{y(x+y)^{-1}}^{(xz+y)(y+x)^{-1}} v^{n-1}\log(1-v)dv.
\end{align}

Therefore, applying \eqref{equ-log-m-1-x}, by a direct calculation, we deduce
\begin{align}
\sum_{k=0}^n x^ky^{n-k}\binom{n}{k} \left(\sum_{j=1}^k \frac{z^j}{j}\right)=(x+y)^n \sum_{j=1}^n \frac{\left(\frac{xz+y}{y+x}\right)^j-\left(\frac{y}{y+x}\right)^j}{j}.
\end{align}
This completes the proof of Theorem \ref{thm:Mneimneh-type-G} (Noting that from $x/(x+y)>0$ and $z\in (-\infty,1)$ gives $(xz+y)(x+y)^{-1}<1$).

Multiplying \eqref{equ:Mneimneh-type-G} by $1/z$ and  integrating over the interval (0,1), we have
\begin{align}\label{cor1-5-one}
\sum_{k=0}^n x^k y^{n-k} \binom{n}{k}H_k^{(2)}&=\sum_{j=1}^n \frac{(x+y)^{n-j}}{j} \int_0^1 \frac{(xz+y)^j-y^j}{z}dz\nonumber\\
&=\sum_{j=1}^n \frac{(x+y)^{n-j}}{j}  y^j \sum_{i=1}^j \int_0^1 \Big(\frac{xz}{y}+1\Big)^{i-1}d\frac{xz}{y}\nonumber\\
&=(x+y)^n \left\{\sum_{j=1}^n \frac{y^j(x+y)^{-j}}{j}\sum_{i=1}^j \frac{y^{-i}(x+y)^i}{i}-\sum_{j=1}^n \frac{y^j(x+y)^{-j}}{j}H_j\right\}.
\end{align}
Then, applying the well-known identity
\begin{align*}
\sum_{j=1}^{n} \frac{H_j}{j}=\zeta^\star_n(1,1)=\frac1{2}Y_2(n)=\frac{H_n^2+H_n^{(2)}}{2}
\end{align*}
and substituting \eqref{cor1-5-one} into \eqref{TMG-Mneimneh-identity} yields \eqref{TMG-Mneimneh-identity-two}. Thus, we complete the proof of Corollary \ref{cor:Mneimneh-type-order2}.

\section{Conclusion}

We presented the explicit formulas of the following Mneimneh-type binomial sum of (alternating) harmonic numbers
\begin{align*}
&\sum_{k=0}^n x^ky^{n-k}\binom{n}{k} Y_p(k),\quad \sum_{k=0}^n x^ky^{n-k}\binom{n}{k} {\bar H}_k,\quad \sum_{k=0}^n x^k y^{n-k} \binom{n}{k}H_k^{(2)}\quad\text{and}\quad \sum_{k=0}^n x^k y^{n-k} \binom{n}{k}H^2_k.
\end{align*}
It is possible that some of other Mneimneh-type binomial sums can be obtained using techniques of the present paper. For example, multiplying \eqref{equ:Mneimneh-type-G} by $\log^{p-1}(z)/z\ (p\in \N)$ and integrating over the interval (0,1), we have
\begin{align}\label{equ-noexp}
\sum_{k=0}^n x^k y^{n-k} \binom{n}{k}H_k^{(p+1)}&=\frac{(-1)^{p-1}}{(p-1)!}\sum_{j=1}^n \frac{(x+y)^{n-j}}{j} \int_0^1 \log^{p-1}(z)\frac{(xz+y)^j-y^j}{z}dz.
\end{align}
Hence, if the evaluation of the integral on the right hand of above can be established, we can obtain the explicit formula of Mneimneh-type binomial sums on the left hand of above. We leave the detail to the interested reader. It should be emphasized that Komatsu-Wang \cite[Eq. (4)]{KW2024} gave an evaluation of \eqref{equ-noexp} with $x=1-q$ and $y=q$, where $q$ is a real number. Gen$\check{\rm c}$ev \cite[Thm. 2.1]{G2024} established a generalization of Mneimneh summation formula
\begin{align}\label{equ-exp-z}
\sum_{k=0}^n \binom{n}{k}\left(\sum_{j=1}^k \frac{z^j}{j^r}\right)p^k(1-p)^k=\sum_{n\geq n_1 \geq \cdots n_r\geq 1} \frac{(1-p)^{n_1}\left(\left(1+\frac{zp}{1-p}\right)^{n_r}-1\right)}{n_1\cdots n_r},
\end{align}
where $r\in \N,p,z\in \R$ and $p\neq 1$.
\medskip

{\bf Declaration of competing interest.}
The authors declares that they has no known competing financial interests or personal relationships that could have
appeared to influence the work reported in this paper.

{\bf Data availability.}
No data was used for the research described in the article.

{\bf Acknowledgments.}  Ce Xu expresses his deep gratitude to Profs. J. Zhao and J.M. Campbell for valuable discussions and comments. Ce Xu is supported by the National Natural Science Foundation of China (Grant No. 12101008), the Natural Science Foundation of Anhui Province (Grant No. 2108085QA01) and the University Natural Science Research Project of Anhui Province (Grant No. KJ2020A0057).


\begin{thebibliography}{99}

\bibitem{C2023}
J.M. Campbell, On Mneimneh's binomial sum involving harmonic numbers, \emph{Discrete Math.} \textbf{346}(2023), 113549.

\bibitem{C1974}
L. Comtet, \emph{Advanced Combinatorics}, D. Reidel Publishing Co., Dordrecht, 1974.

\bibitem{CG1996}
J.H. Conway and R.K. Guy, \emph{The book of numbers}. Springer, New York (1996).

\bibitem{G2024}
M. Gen$\check{\rm c}$ev, Generalized Mneimneh sums and their application to multiple
polylogarithms via Toeplitz theorem, private email dated March 9, 2024. 

\bibitem{H1992}
M.E. Hoffman, Multiple harmonic series, \emph{Pacific J. Math.} \textbf{152}(1992), 275--290.

\bibitem{KW2024}
T. Komatsu and P. Wang, A generalization of Mneimneh's binomial sum of harmonic numbers, \emph{Discrete Math.} \textbf{347}(2024), 113549.

\bibitem{M2023}
S. Mneimneh, A binomial sum of harmonic numbers, \emph{Discrete Math.} \textbf{346}(2023), 113075.

\bibitem{Riordan58}
J. Riordan, \emph{An Introduction to Combinatorial Analysis}, Reprint of the 1958 original, Dover Publications, Inc., Mineola, NY, 2002.

\bibitem{Xu2017}
C. Xu, Multiple zeta values and Euler sums, \emph{J. Number Theory} \textbf{177}(2017), 443--478.

\bibitem{XYS2016}
C. Xu, Y. Yan and Z. Shi, Euler sums and integrals of polylogarithm functions, \emph{J. Number Theory} \textbf{165}(2016), 84-108.

\bibitem{XuZhao2020d}
C.\ Xu and J. Zhao, Explicit relations of some variants of convoluted multiple zeta values. arXiv:2103.01377.

\bibitem{DZ1994}
D. Zagier, \emph{Values of zeta functions and their applications}, in: First European Congress
of Mathematics, Vol.~II, pp.\ 497--512, Birkhauser, Boston, 1994.

\end{thebibliography}
\end{document}